# Low-Thrust Transfer Between Circular Orbits Using Natural Precession

Max CERF [*]

Airbus Defence and Space, France

## Abstract

The minimum-fuel low-thrust transfer between circular orbits is formulated using the Edelbaum's averaged dynamics with the addition of the nodal precession due to the first zonal term. The extremal analysis shows that an optimal transfer is composed of three sequences in the regular case. The optimal control problem is solved by a shooting method with a costate guess derived from an approximate solution.

## 1.   Introduction

The low-thrust technology for orbit transfers offers promising performance in terms of fuel consumption. This fuel gain is nevertheless offset by a longer transfer duration that may be penalizing for commercial Earth Orbit Raising (EOR) missions, for example to raise a telecommunication satellite from a Low Earth Orbit (LEO) to a Geostationary Earth Orbit (GEO). A particular application where time is not an issue is the removal of spent satellites from the near Earth region. Active Debris Removal (ADR) missions aiming at deorbiting several debris are currently under study. For such missions the low-thrust technology proves particularly attractive and a quick performance assessment method would be very useful at the vehicle design stage. The targeted debris are mostly old observation satellites evolving on near circular Sun-Synchronous Orbits. The transfer strategy between two debris must account not only for the altitude and inclination change, but also for the ascending node precession due to the Earth flattening (first zonal term $J_2$).

Several analytical or quasi-analytical solutions exist for the minimum-time low-thrust transfer problem between circular orbits. Some of them rely on perturbation techniques allowing an explicit solution under specific assumptions [1]. The most famous solution is due to Edelbaum in 1961 [2,3]. It is based on an averaged dynamical model assuming a constantly circular orbit and continuous thrusting. The original work of Edelbaum has been retrieved in the frame of optimal control theory [4] and several extensions have been derived in order to enlarge the application scope of the model.

- In [5] the Edelbaum's model is enhanced to account for the mass variation during the transfer and for a variable specific impulse allowing a reduction of the fuel consumption.

---

[*] AIRBUS Defence & Space, 78130 Les Mureaux, France
max.cerf@astrium.eads.net

- In [6] the solution of the Edelbaum's problem is enhanced to comply with an altitude upper bound. This is particularly useful when a large inclination change is required since Edelbaum's solution involves high altitudes to minimize the cost of the inclination change.
- Still in [6] the dynamics averaging method is adapted to perform a RAAN change instead of an inclination change.
- In [7] no-thrust legs due to Earth shadowing are accounted by shifting accordingly the trajectory time at each revolution.

These extensions are devoted to the minimum-time problem. For the debris removal missions under interest, the goal is rather to minimize the fuel consumption within a fixed duration allocated to the mission, assuming that this duration is sufficiently large wrt the low-thrust engine capabilities.

This paper addresses the minimum-fuel low-thrust transfer problem between circular orbits taking into account a Right Ascension Ascending Node (RAAN) constraint. The optimal control problem is formulated using Edelbaum's dynamical model and the solution structure is analyzed. A solution method is proposed using an indirect method and a costate guess derived from a simplified problem. The method is exemplified on an application case representative of a debris removal mission.

## 2. Problem Formulation and Analysis

This section formulates the low-thrust transfer problem in terms of an Optimal Control Problem (OCP). The solution structure is then analyzed by applying the Pontryaguin Maximum Principle (PMP).

### 2.1 Dynamics

The dynamics is based on the Edelbaum's dynamical model which assumes constantly circular orbits. At a given date the orbital plane is defined in the Earth inertial reference frame by the inclination I and the right ascension of the ascending node $\Omega$ (RAAN). The inclination I is the angle of the orbital plane with the Earth equatorial plane. The intersection of the orbital plane with the Equator is the line of nodes. The RAAN $\Omega$ is the angle between the X axis of the Earth inertial reference frame and the direction of the ascending node (node crossed with a northwards motion). The circular orbit shape is defined equivalently by the radius a or by the velocity V.

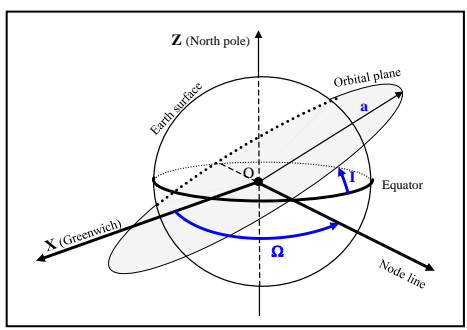

Figure 1 : Orbital parameters for a circular orbit

a : radius
V : velocity
I : inclination
$\Omega$ : right ascension of the ascending node (RAAN)

The Edelbaum's model is derived by an averaging of the Gauss equations, with the following assumptions :
- The evolution of the orbital parameters is averaged on one period.
- The averaged orbit is constantly circular throughout the transfer.
- The acceleration level denoted f is constant throughout the transfer.
- The thrust direction is normal to the radius vector and it makes a constant angle β with the orbital plane during one period with a sign change at the antinodes. This results in a null RAAN change over the period while maximizing the inclination change for the current value of the out of plane angle β.

In this paper the upper case notations V, I, Ω represent the averaged orbital parameters.

The Edelbaum's problem consists in minimizing the transfer duration. The engine is continuously thrusting with the constant acceleration level f and the transfer is controlled by varying the averaged angle β from one period to the other. We refer to [4] and [8] for the detailed model formulation and the analytical solution of the minimum-time problem.

For our purpose the Edelbaum's dynamics is extended by taking into account the first zonal term (denoted $J_2$) due to the Earth flattening. This perturbation causes no secular change on the semi-major axis, the eccentricity and the inclination. The averaged orbit remains circular with radius a and inclination I. On the other hand, a RAAN precession rate is induced depending on the orbit radius a and the inclination I [8].

$$\dot{\Omega} = -\frac{3}{2} J_2 \sqrt{\mu} R_E^2 \, a^{-\frac{7}{2}} \cos I = -k \left(\frac{\mu}{a}\right)^{\frac{7}{2}} \cos I = -kV^7 \cos I \tag{1}$$

The constants of the Earth gravitational model are [8] :  $R_E = 6378137$ m   (equatorial radius)

$\mu = 3.986005 \cdot 10^{14}$ m$^3$/s$^2$   (gravitational constant)

$J_2 = 1.08266$   (first zonal term)

$$k = \frac{3 J_2 R_E^2}{2 \mu^3} = 1.0425 \times 10^{-33} \text{ s}^6 / \text{m}^7$$

The dynamics model consists thus in three ordinary differential equations representing the evolution of the averaged orbital parameters (V, I, Ω).

$$\begin{cases} \dot{V} = -f \cos\beta \\ \dot{I} = \dfrac{2}{\pi V} f \sin\beta \\ \dot{\Omega} = -kV^7 \cos I \end{cases} \tag{2}$$

The control variables are the acceleration level f and its direction β over each period. The acceleration level can be varied between zero and a maximum value $f_{max}$.

## 2.2   Optimal Control Problem

The problem we are interested in consists in transferring the vehicle from a given circular orbit to another given circular orbit in a given duration, while minimizing the fuel consumption.

The fuel consumed is linked to the velocity impulse through the rocket equation derived by Tsiolkovsky [8].

$$\Delta V = \int_{t_0}^{t_f} f \, dt = v_e \ln \frac{M_0}{M_f} \quad \Rightarrow \quad m_c = M_0 - M_f = M_0 \left(1 - e^{-\frac{\Delta V}{v_e}}\right) \tag{3}$$

$M_0$ and $M_f$ are respectively the initial and final gross mass, $v_e$ is the engine exhaust velocity and $m_c$ is the propellant mass consumed. Minimizing the fuel consumption $m_c$ is equivalent to minimizing the velocity impulse $\Delta V$. The Optimal Control Problem (OCP) is formulated considering the velocity impulse as cost function, with completely prescribed endpoints and bounded control.

Optimal Control Problem

$$\underset{f,\beta}{\text{Min}} \; J = \int_{t_0}^{t_f} f \, dt \quad \text{s.t.} \quad \begin{cases} \dot{V} = -f \cos\beta \\ \dot{I} = \dfrac{2}{\pi V} f \sin\beta \\ \dot{\Omega} = -k V^7 \cos I \end{cases} \tag{4}$$

Fixed initial conditions : $t_0$, $V(t_0)=V_0$, $I(t_0)=I_0$, $\Omega(t_0)=\Omega_0$

Fixed final conditions : $t_f$, $V(t_f)=V_f$, $I(t_f)=I_f$, $\Omega(t_f)=\Omega_f$

Control bounds : $0 \leq f \leq f_{max}$

The problem is autonomous since the dynamics and the integral cost do not depend explicitly on the time. On the other hand the endpoint conditions $\Omega(t_0)$ and $\Omega(t_f)$ are time-dependent due to the natural precession. It can be observed that the velocity and the inclination evolutions do not depend on the RAAN value. The problem is insensitive to a RAAN shift applied identically on $\Omega(t_0)$ and $\Omega(t_f)$.

The thrust is used to modify directly the velocity and the inclination whilst the RAAN change is performed passively by the natural $J_2$ precession. This control strategy based on Edelbaum's model is sub-optimal due to the dynamics simplifications. It should be envisioned only if a sufficiently large duration is allocated to the mission allowing benefiting from the $J_2$ precession. This is normally the case for practical applications such as debris removal missions [9]. If not, the Edelbaum's control model is no longer suited and the problem should be addressed without control law simplifications.

### 2.3 Extremal Analysis

According to the Pontryaguin Maximum Principle (PMP) [10,11], the necessary conditions for $(f,\beta)$ to be an optimal control for the OCP are the existence of absolutely continuous functions $p_V$, $p_I$, $p_\Omega$ (respective costates of V, I, $\Omega$), and the existence of a non positive real $p_0$ (cost multiplier) such that :

- The control $(f, \beta)$ maximizes nearly everywhere the Hamiltonian function defined as :

$$H = p_0 f + p_V \dot{V} + p_I \dot{I} + p_\Omega \dot{\Omega} \tag{5}$$

- The costate vector $(p_V, p_I, p_\Omega)$ is not identically null on an interval of $[t_0,t_f]$ and it satisfies the differential system :

$$\begin{cases} \dot{p}_V = -\dfrac{\partial H}{\partial V} = 7p_\Omega kV^6 \cos I + \dfrac{2}{\pi V^2} p_I f \sin\beta \\ \dot{p}_I = -\dfrac{\partial H}{\partial I} = -p_\Omega kV^7 \sin I \\ \dot{p}_\Omega = -\dfrac{\partial H}{\partial \Omega} = 0 \end{cases} \qquad (6)$$

Since the initial and final conditions are completely fixed, there are no transversality conditions on the costate and on the Hamiltonian. Introducing the switching function S :

$$S = p_0 - p_V \cos\beta + p_I \dfrac{2}{\pi V} \sin\beta \qquad (7)$$

the Hamiltonian is expressed as :

$$H = fS + p_\Omega \dot{\Omega} \qquad (8)$$

The term $p_\Omega \dot{\Omega}$ does not depend explicitly on the control variables (f,β). The Hamiltonian is linear wrt the acceleration level f and it depends on the thrust direction β through the switching function S.

The Hamiltonian maximization yields the optimal acceleration level depending on the switching function sign.

When S > 0 then f = $f_{max}$    (thrust arc)

When S < 0 then f = 0        (coast arc)    (9)

When S = 0, f cannot be determined directly from the PMP (singular arc)

In the case of a non null acceleration, the Hamiltonian maximization wrt the thrust direction β yields :

$$\begin{cases} \dfrac{\partial H}{\partial \beta} = 0 \quad \Rightarrow \quad p_V \sin\beta + p_I \dfrac{2}{\pi V} \cos\beta = 0 \\ \dfrac{\partial^2 H}{\partial \beta^2} \leq 0 \quad \Rightarrow \quad p_V \cos\beta - p_I \dfrac{2}{\pi V} \sin\beta \leq 0 \end{cases} \qquad (10)$$

The costates $p_V$ and $p_I$ can not vanish simultaneously on a time interval, else from Eq. (6) : $p_\Omega$=0 which would be a contradiction to the PMP (non trival costate). The system Eq. (10) defines therefore the thrust direction β without sign ambiguity. Using for example the function Atan2(s,c) where s and c are respectively proportional to the angle sine and cosine, we can express β as :

$$\beta = \text{Atan2}\left(\dfrac{2}{\pi V} p_I, -p_V\right) \qquad (11)$$

Before investigating the solution structure, several preliminary results are derived here after.

Property 1 :  Abnormal solution ($p_0$=0)
The inequality in Eq. (10) gives a lower bound on the switching function Eq. (7).

$$S = p_0 - p_V \cos\beta + p_I \dfrac{2}{\pi V} \sin\beta \geq p_0 \qquad (12)$$

The solution is said abnormal when the cost multiplier $p_0$ is null. The switching function is then positive or null and the trajectory is made of maximum thrust arcs and singular arcs. The singular case is studied further in §2.5.

Property 2 : Value of H, $p_\Omega$ and fS

The OCP Eq. (4) is autonomous so that the Hamiltonian is constant. From Eq. (6), $p_\Omega$ is also constant and from Eq. (9) the term fS is positive or null.

$$\begin{cases} H(t) = fS + p_\Omega \dot{\Omega} = C^{te} \\ p_\Omega = C^{te} \\ fS \geq 0 \end{cases} \qquad (13)$$

Property 3 : Case $p_\Omega \neq 0$

If $p_\Omega \neq 0$, we define $\dot{\Omega}_d$ by $H \underset{def}{=} p_\Omega \dot{\Omega}_d$. H and $p_\Omega$ being constant from Eq. (13), $\dot{\Omega}_d$ is also constant.

Using $fS \geq 0$ from Eq. (13), we have :

$$p_\Omega \dot{\Omega} = p_\Omega \dot{\Omega}_d - fS \leq p_\Omega \dot{\Omega}_d \qquad (14)$$

The term $p_\Omega \dot{\Omega}$ is bounded upperly by $p_\Omega \dot{\Omega}_d$ and the bound is reached whenever fS = 0.

All coast arcs (f=0) and all singular arcs (S=0) have therefore the same precession rate equal to $\dot{\Omega}_d$.

Property 4 : Case $p_\Omega = 0$

If $p_\Omega=0$, we have $H = fS \geq 0$ from Eq. (13). We consider successively the case $H \neq 0$ and the case H=0.
- If H > 0, then S > 0 and from Eq. (9) the acceleration level is constantly maximum : $f = f_{max}$.
  The OCP Eq. (4) is equivalent to Edelbaum's minimum-time problem whose solution is analytical [4,8].
- If H = 0, then either f=0 (coast arc), or S=0 (singular arc). The singular case is studied further in §2.5.

A regular solution with $p_\Omega=0$ is therefore composed of a single arc either at f=0 or $f=f_{max}$.

Property 5 : Cost derivatives wrt initial and final date

If the initial date is shifted from $t_0$ to $t_0'$, the velocity and inclination are unchanged whereas the RAAN is modified by the natural precession.

$$\Omega(t_0') = \Omega(t_0) + \dot{\Omega}(t_0).(t_0' - t_0) \quad \text{with} \quad \dot{\Omega}(t_0) = -kV(t_0)^7 \cos I(t_0) \qquad (15)$$

The total derivative of the optimal cost denoted J* wrt the initial date is :

$$\frac{dJ^*}{dt_0} = \frac{\partial J^*}{\partial t_0} + \frac{\partial J^*}{\partial \Omega(t_0)} \cdot \frac{d\Omega(t_0)}{dt_0} = \frac{\partial J^*}{\partial t_0} + \frac{\partial J^*}{\partial \Omega(t_0)} \cdot \dot{\Omega}(t_0) \qquad (16)$$

The partial derivatives of the optimal cost wrt to the initial date and the RAAN are linked respectively to the Hamiltonian and to the RAAN costate [11].

$$\frac{\partial J^*}{\partial t_0} = H(t_0) \quad , \quad \frac{\partial J^*}{\partial \Omega(t_0)} = -p_\Omega(t_0) \qquad (17)$$

Replacing in Eq. (16) with H given by Eq. (8), we get after simplification :

$$\frac{dJ^*}{dt_0} = f(t_0)S(t_0) \geq 0 \tag{18}$$

Similarly for a final date shift from $t_f$ to $t_f$', we have with opposite signs for the partial derivatives :

$$\frac{\partial J^*}{\partial t_f} = -H(t_f) \quad , \quad \frac{\partial J^*}{\partial \Omega(t_f)} = p_\Omega(t_f) \tag{19}$$

$$\frac{dJ^*}{dt_f} = -f(t_f)S(t_f) \leq 0 \tag{20}$$

Eq. (18) has two physical interpretations.

- Assume that an optimal trajectory denoted T*($t_0$) has been found for OCP with the cost value J*($t_0$). We consider the problem OCP' identical to OCP but starting at an earlier initial date $t_0$' < $t_0$ with $\Omega(t_0$') given by Eq. (15). A feasible trajectory for OCP' consists in a zero thrust sequence from $t_0$' to $t_0$ followed by the trajectory T*($t_0$) from $t_0$ to $t_f$. This feasible trajectory has the cost J*($t_0$), so that the optimal cost J*($t_0$') for OCP' is at most equal to J*($t_0$). The cost total derivative wrt the initial date is therefore positive.

$$\begin{cases} t_0' < t_0 \\ J^*(t_0') \leq J^*(t_0) \end{cases} \Rightarrow \frac{dJ^*}{dt_0} \geq 0 \tag{21}$$

- Assume that the optimal trajectory starts with a coast arc from $t_0$ to $t_1$. From Eq. (18) with f($t_0$)=0, the optimal cost is insensitive to the initial date. Restating OCP Eq. (4) with any initial date $t_0$' comprised between $t_0$ and $t_1$ yields an equivalent problem having the same optimal cost J*($t_0$) as OCP.

The interpretations are similar for Eq. (20) at the final date.

The control structure is next investigated successively in the regular case and in the singular case.

## 2.4 Regular Solution

In the regular case, the switching function S does not vanish identically on any interval of [$t_0$;$t_f$]. The acceleration level is either 0 or $f_{max}$ depending on the sign of S. If $p_\Omega$=0 the solution is composed of a single arc (Property 4) and the problem has a direct analytical solution. Assuming $p_\Omega \neq 0$ in the sequel of this section, all coast phases have the same precession rate $\dot{\Omega}_d$ (Property 3). During the coast phases the velocity and the inclination remain constant, whereas the RAAN evolves linearly with the precession rate $\dot{\Omega}_d$.

The Figure 2 illustrates a possible evolution of $p_\Omega \dot{\Omega}$ along the transfer, with 3 propelled sequences separated by 2 coast phases. This scenario is called scenario 1.

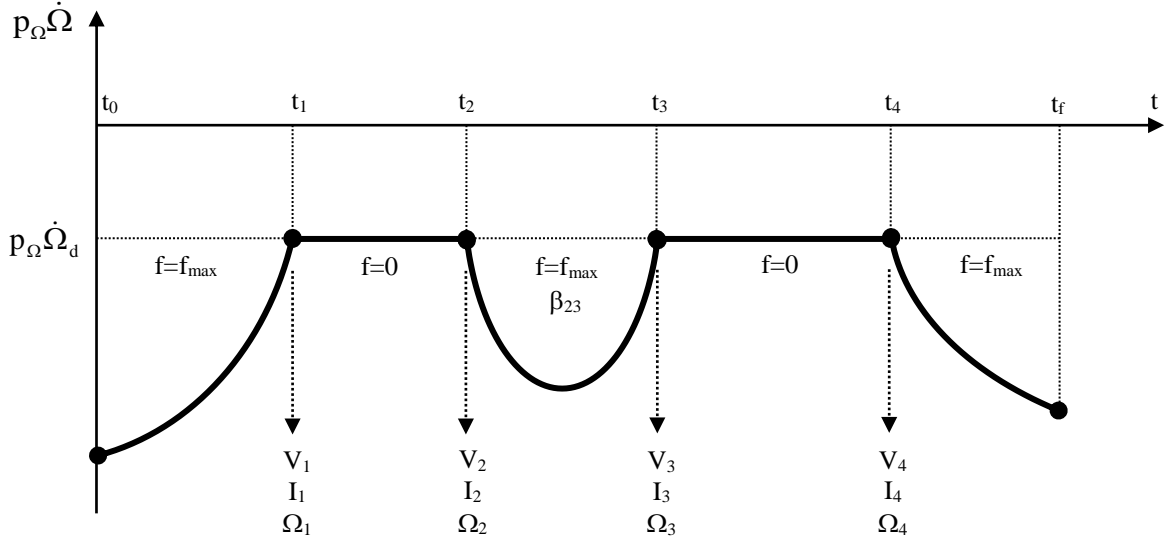

Figure 2 : Scenario 1

The following relationships hold between the endpoints of the coast phases.

$$\begin{cases} V_2 = V_1 \\ I_2 = I_1 \\ \Omega_2 = \Omega_1 + \dot{\Omega}_d \cdot (t_2 - t_1) \end{cases} \quad \text{and} \quad \begin{cases} V_4 = V_3 \\ I_4 = I_3 \\ \Omega_4 = \Omega_3 + \dot{\Omega}_d \cdot (t_4 - t_3) \end{cases} \quad (22)$$

The leg from $t_2$ to $t_3$ is at maximum acceleration level $f_{max}$ with a thrust direction law denoted $\beta_{23}$. This control transfers the state from $(V_2,I_2,\Omega_2)$ to $(V_3,I_3,\Omega_3)$ in a duration $\Delta t_{23} = t_3 - t_2$. The RAAN change during this leg is denoted $\Delta\Omega_{23} = \Omega_3 - \Omega_2$.

We consider now the scenario obtained by permuting the 3rd and the 4th leg, each leg keeping the same duration and the same control law as in the scenario 1. The date $t_2$ and $t_4$ are unchanged, whilst the intermediate date $t_3$ is shifted to $t_3$'. The new scenario called scenario 2 is depicted on the Figure 3.

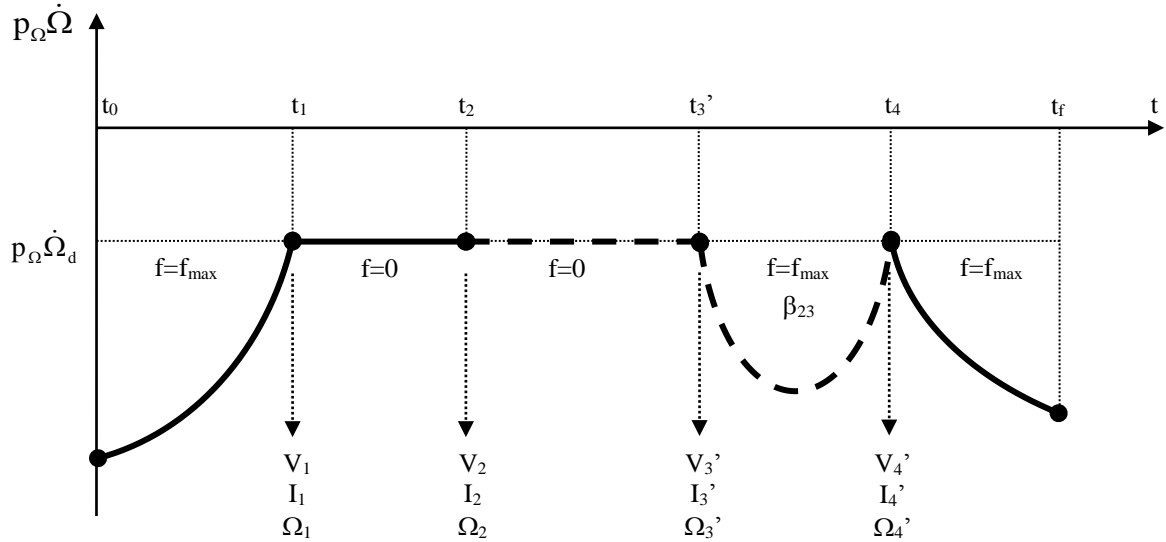

Figure 3 : Scenario 2

The scenario 2 is identical to the scenario 1 until the date $t_2$. The coast phase starting at $t_2$ has the duration $\Delta t_{34}=t_4-t_3$ and the precession rate $\dot{\Omega}_d$. It ends at $t_3$' with the state :

$$\begin{cases} V_3' = V_2 \\ I_3' = I_2 \\ \Omega_3' = \Omega_2 + \dot{\Omega}_d.\Delta t_{34} \end{cases} \quad (23)$$

The control ($f_{max}$, $\beta_{23}$) retrieved from the scenario 1 is applied from $t_3$' to $t_4$. Compared to the leg from $t_2$ to $t_3$ of the scenario 1, this leg from $t_3$' to $t_4$ has a different starting date ($t_3$' instead of $t_2$), the same duration ($\Delta t_{23}$) and the same initial state ($V_3$'=$V_2$, $I_3$'=$I_2$) except for the RAAN ($\Omega_3$' instead of $\Omega_2$).

The shift of the starting date does not change the trajectory since the problem is autonomous. As observed at the end of §2.2, the solution is also insensitive to a RAAN shift. Applying the control ($f_{max}$, $\beta_{23}$) from the initial state ($V_2,I_2,\Omega_3$') at the date $t_3$' until the date $t_4$ yields a leg identical to scenario 1 with a RAAN change equal to $\Delta\Omega_{23}$. The state at $t_4$ is thus :

$$\begin{cases} V_4' = V_3 \\ I_4' = I_3 \\ \Omega_4' = \Omega_3' + \Delta\Omega_{23} = \Omega_2 + \dot{\Omega}_d.\Delta t_{34} + \Delta\Omega_{23} = \Omega_3 + \dot{\Omega}_d.\Delta t_{34} \end{cases} \quad (24)$$

Comparing with Eq. (22), the final state at $t_4$ is identical to scenario 1.

$$\begin{cases} V_4' = V_4 \\ I_4' = I_4 \\ \Omega_4' = \Omega_4 \end{cases} \quad (25)$$

The scenario 2 reaches the same final conditions as the scenario 1. Their costs are identical since the thrusting sequences at $f_{max}$ have the same duration. A solution with several coast phases can thus be replaced by an equivalent solution with a single coast phase, and we have the following result.

Property 6 : Structure of a regular solution
In the regular case, a solution of the optimal control problem can be sought assuming a 3 sequences control structure of the form : $f_{max} - 0 - f_{max}$.

### 2.5 Singular Solution

In the singular case, the switching function S vanishes identically on an interval of [$t_0$;$t_f$]. The acceleration level can no longer be deduced directly from the Hamiltonian maximization.
Combining Eq. (7) with Eq. (10) yields the costates $p_V$ and $p_I$ :

$$\begin{cases} S=0 \\ \dfrac{\partial H}{\partial \beta}=0 \end{cases} \Rightarrow \begin{cases} p_0 - p_V \cos\beta + p_I \dfrac{2}{\pi V}\sin\beta = 0 \\ p_V \sin\beta + p_I \dfrac{2}{\pi V}\cos\beta = 0 \end{cases} \Rightarrow \begin{cases} p_0 \cos\beta - p_V = 0 \\ p_0 \sin\beta + p_I \dfrac{2}{\pi V} = 0 \end{cases} \Rightarrow \begin{cases} p_V = p_0 \cos\beta \\ p_I = -\dfrac{\pi V}{2} p_0 \sin\beta \end{cases} \quad (26)$$

It was stated (Property 1) that an abnormal solution ($p_0=0$) comprises either maximum thrust arcs or singular arcs. Assuming $p_0=0$ in Eq. (26) leads to $p_V=0$, $p_I=0$ and using Eq. (6) to $p_\Omega=0$. The costates vanish identically which is in contradiction to the PMP. A singular solution cannot be abnormal and we have the following result.

Property 7 : Abnormal solution

An abnormal solution is composed of a single maximum thrust arc. In that case, the OCP Eq. (4) becomes equivalent to Edelbaum's minimum-time problem whose solution is analytical [4,8].

The usual way to analyze a singular solution consists in differentiating the switching function Eq. (7) :

$$\dot{S} = -\dot{p}_V \cos\beta + p_V \dot\beta \sin\beta + \dot{p}_I \frac{2}{\pi V}\sin\beta - p_I \frac{2\dot V}{\pi V^2}\sin\beta + p_I \frac{2}{\pi V}\dot\beta \cos\beta$$
$$= -\dot{p}_V \cos\beta + \dot{p}_I \frac{2}{\pi V}\sin\beta - p_I \frac{2\dot V}{\pi V^2}\sin\beta + \dot\beta\left(p_V \sin\beta + p_I \frac{2}{\pi V}\cos\beta\right) \quad (27)$$

The parenthesis is null from Eq. (10). Using Eq. (4) and Eq. (6) to replace the derivatives of V, $p_V$, $p_I$ :

$$\dot{S} = -\left(7 p_\Omega k V^6 \cos I + \frac{2}{\pi V^2} p_I f \sin\beta\right)\cos\beta - p_\Omega k V^7 \sin I \frac{2}{\pi V}\sin\beta + p_I \frac{2}{\pi V^2} f \cos\beta \sin\beta \quad (28)$$

After simplification we get :

$$\dot{S} = -p_\Omega k V^6\left(7\cos I \cos\beta + \frac{2}{\pi}\sin I \sin\beta\right) \quad (29)$$

Along a singular arc, the switching function derivative vanishes identically.

$$\dot{S}=0 \Rightarrow p_\Omega=0 \quad \text{or} \quad 7\cos I\cos\beta + \frac{2}{\pi}\sin I\sin\beta = 0 \quad (30)$$

We consider successively the case $p_\Omega=0$ and $p_\Omega\neq 0$.

Case $p_\Omega=0$

If $p_\Omega=0$, the derivative of S is null whatever the date. The trajectory consists in a single singular arc spanning from $t_0$ to $t_f$. The costate equations Eq. (6) become identical to the Edelbaum's minimum-time problem, with the difference that the acceleration level is no longer a constant. Edelbaum's formulae can be partly retrieved following a similar approach to [4,8].

From Eq. (6), the costate $p_I$ is constant, which implies from Eq. (26) that the product $V\sin\beta$ is constant along the singular arc.

$$V\sin\beta = V_0 \sin\beta_0 \quad (31)$$

Differentiating Eq. (31) and using Eq. (4) to replace $\dot V$ yields an expression for $\dot\beta$ :

$$\dot V \sin\beta + V\dot\beta \cos\beta = 0 \Rightarrow \dot\beta = \frac{f\sin\beta}{V} \quad (32)$$

This expression appears in the right hand side of Eq. (4) for $\dot{I}$. Replacing and integrating yields an analytical expression for the inclination I depending on the control β along the singular arc.

$$\dot{I} = \frac{2}{\pi V} f \sin\beta = \frac{2}{\pi}\dot{\beta} \quad \Rightarrow \quad I - I_0 = \frac{2}{\pi}(\beta - \beta_0) \tag{33}$$

Using Eq. (33) to replace β in Eq. (31) yields the value of $\beta_0$ depending on the endpoint velocity and inclination.

$$V\sin\beta = V\sin\left(\beta_0 + \frac{\pi}{2}(I - I_0)\right) = V_0 \sin\beta_0 \quad \Rightarrow \quad \tan\beta_0 = \frac{V_f \sin\left(\frac{\pi}{2}(I - I_0)\right)}{V_0 - V_f \cos\left(\frac{\pi}{2}(I - I_0)\right)} \tag{34}$$

The cost of the singular solution can also be assessed analytically starting from the velocity equation Eq. (4).

$$\frac{dV}{dt} = -f\cos\beta \quad \Rightarrow \quad fdt = -\frac{dV}{\cos\beta} = \pm\frac{dV}{\sqrt{1-\sin^2\beta}} = \pm\frac{VdV}{\sqrt{V^2 - V^2\sin^2\beta}} = \pm\frac{VdV}{\sqrt{V^2 - V_0^2 \sin^2\beta_0}} \tag{35}$$

Integrating from $t_0$ to $t_f$ yields the cost value.

$$J = \int_{t_0}^{t_f} fdt = \pm\int_{t_0}^{t_f} \frac{VdV}{\sqrt{V^2 - V_0^2 \sin^2\beta_0}} = \pm\left[\sqrt{V^2 - V_0^2 \sin^2\beta_0}\right]_{t_0}^{t_f} = \pm\sqrt{V_f^2 - V_0^2 \sin^2\beta_0} \pm V_0 \cos\beta_0 \tag{36}$$

The cost value only depends on the endpoint conditions and not on the acceleration level along the trajectory.
In order to completely define the singular trajectory there remains to find the variation of acceleration level f(t) from $t_0$ to $t_f$. The function f(t) defines the thrust direction β(t) by integrating Eq. (32) which in turn defines the velocity Eq. (31), the inclination Eq. (33) and the RAAN by integrating Eq. (1).

No explicit solution has been found for the expression of the acceleration level, so that a direct method seems necessary considering a discretization of the unknown function f(t).
Nevertheless a simplifying approach is possible by assuming a constant acceleration level. Indeed the cost value is insensitive to the function f(t) from Eq. (36). The assumption of a constant acceleration level allows applying the Edelbaum's analytical formulae. An intermediate coast phase is inserted within the Edelbaum's trajectory so that the targeted final conditions can be met at the prescribed final date $t_f$. This coast phase allows controlling the RAAN final value without changing the velocity and the inclination. The problem reduces to finding the acceleration level and the coast phase dates that yields the targeted final conditions at the prescribed final date $t_f$.

Case $p_\Omega \neq 0$
If $p_\Omega \neq 0$, we have from Eq. (30):

$$7\cos I \cos\beta + \frac{2}{\pi}\sin I \sin\beta = 0 \tag{37}$$

Before further analyzing the singular solution, we consider two special cases.
- If $\cos\beta = 0$, then from Eq. (37) we have $\sin I = 0$ so that I = 0 or 180 deg. On the other hand $\sin\beta = \pm 1$, and from Eq. (4) $\dot{I} \neq 0$. The inclination can therefore not remain constantly equal to 0 or 180 deg on the singular arc. From Eq. (37) this means also that $\cos\beta$ cannot remain constantly equal to 0. This can occur only at isolated dates on the singular arc.

- If $\cos I = 0$, then from Eq. (37) we have $\sin\beta = 0$ and from Eq. (4) $\dot{I} = 0$. The inclination is therefore constantly equal to 90 deg. The precession rate is null from Eq. (1) meaning that the RAAN is constant. The problem reduces to a planar transfer between polar circular orbits.

We assume in the sequel that $\cos\beta \neq 0$ and $\cos I \neq 0$. Under these assumptions, Eq. (37) can be rewritten as :

$$\tan I \tan \beta = -\frac{7\pi}{2} \qquad (38)$$

Differentiating Eq. (38) and using Eq. (4) to replace $\dot{I}$ yields an expression for $\dot{\beta}$ :

$$\frac{\dot{I}}{\cos^2 I}\tan\beta + \tan I \frac{\dot{\beta}}{\cos^2 \beta} = 0 \;\;\Rightarrow\;\; \dot{\beta} = -\frac{2}{\pi V} f \frac{\sin^2\beta \cos\beta}{\sin I \cos I} \qquad (39)$$

Differentiating the costate $p_V$ from Eq. (26) and identifying to Eq. (6) :

$$\dot{p}_V = -p_0 \dot{\beta}\sin\beta = 7 p_\Omega k V^6 \cos I + \frac{2}{\pi V^2} p_I f \sin\beta \qquad (40)$$

Replacing $p_I$ with Eq. (26) and $\dot{\beta}$ with Eq. (39) :

$$p_0 f \sin^2 \beta \left(1 + \frac{2}{\pi}\frac{\sin\beta\cos\beta}{\sin I \cos I}\right) = 7 p_\Omega k V^7 \cos I \qquad (41)$$

The right hand side can be simplified using Property 3. Indeed when $p_\Omega \neq 0$, the precession rate is constant along the singular arc : $\dot{\Omega} = -kV^7 \cos I = \dot{\Omega}_d$. We obtain for the acceleration level.

$$f = \frac{-7 p_\Omega \dot{\Omega}_d}{p_0 \sin^2 \beta \left(1 + \frac{2}{\pi}\frac{\sin\beta\cos\beta}{\sin I \cos I}\right)} \qquad (42)$$

Using Eq. (38), $\beta$ can be eliminated in order to assess f as a function of I only. We use the following notations :

$$\begin{cases} \tau = \tan I \\ \alpha = -\frac{7\pi}{2} \end{cases} \qquad (43)$$

The thrust direction Eq. (38) is expressed as :

$$\tan\beta(\tau) = \frac{\alpha}{\tau} \qquad (44)$$

After calculation the acceleration level Eq. (42) along the singular arc is expressed as :

$$f(\tau) = \frac{-p_\Omega \dot{\Omega}_d}{p_0 \alpha^2} \frac{(\tau^2 + \alpha^2)^2}{6\tau^2 - (\alpha^2 - 7)} \qquad (45)$$

Eqs. (44,45) define the singular control depending on the current inclination I and on the parameters $p_0$, $p_\Omega$, $\dot{\Omega}_d$.

The acceleration level Eq. (45) must be bounded between 0 and $f_{max}$. We now analyze the variation of $f(\tau)$ in order to derive conditions on the existence of a singular arc.

The denominator in Eq. (45) vanishes for the value $\tau_s$ corresponding to the inclination values $I_{s1}$ and $I_{s2}$.

$$\tau_s = \pm\sqrt{\frac{\alpha^2 - 7}{6}} \quad \Rightarrow \quad I_{s1} \approx 77.07 \text{ deg or } I_{s2} \approx 102.93 \text{ deg} \tag{46}$$

Differentiating $f(\tau)$ wrt $\tau$, we have :

$$f'(\tau) = \frac{-p_\Omega \dot{\Omega}_d}{p_0 \alpha^2} \frac{\tau(3\tau^2 - 4\alpha^2 + 7)(\tau^2 + \alpha^2)}{(6\tau^2 - (\alpha^2 - 7))^2} \tag{47}$$

The derivative $f'(\tau)$ vanishes either for $\tau=0$ (corresponding to I=0 or 180 deg), or for the value $\tau_m$ corresponding to the inclination values $I_{m1}$ and $I_{m2}$.

$$\tau_m = \pm\sqrt{\frac{4\alpha^2 - 7}{3}} \quad \Rightarrow \quad I_{m1} \approx 85.46 \text{ deg or } I_{m2} \approx 94.54 \text{ deg} \tag{48}$$

The respective values of the acceleration level are :

$$\begin{cases} f_0 = f(\tau = 0) = \dfrac{p_\Omega \dot{\Omega}_d}{p_0} \dfrac{\alpha^2}{\alpha^2 - 7} \\ f_m = f(\tau = \tau_m) = -\dfrac{7}{9} \dfrac{p_\Omega \dot{\Omega}_d}{p_0} \dfrac{\alpha^2 - 1}{\alpha^2} \end{cases} \tag{49}$$

The variation of the acceleration level $f(\tau)$ is depicted in the Table 1 in the case $p_\Omega \dot{\Omega}_d > 0$.

| I (deg) | 0 | | $I_{s1} \approx 77.1$ | | $I_{m1} \approx 85.5$ | | 90 | | $I_{m2} \approx 94.5$ | | $I_{s2} \approx 102.9$ | | 180 |
|---|---|---|---|---|---|---|---|---|---|---|---|---|---|
| f' | | − | | − | | + | | − | | + | | + | |
| f | $f_0 \searrow -\infty$ | | $+\infty \searrow f_m$ | | $f_m \nearrow +\infty$ | | $+\infty \searrow f_m$ | | $f_m \nearrow +\infty$ | | $-\infty \nearrow f_0$ | | |

Table 1 : Variation of the acceleration level in the case $p_\Omega \dot{\Omega}_d > 0$

We continue the analysis in the case $p_\Omega \dot{\Omega}_d > 0$ (the case $p_\Omega \dot{\Omega}_d < 0$ is discussed after).

Since the cost multiplier $p_0$ is negative, we have in Eq. (49) : $f_0 < 0$, $f_m > 0$. The singular arc may exist only in an inclination range around $I_{m1}$ or $I_{m2}$. Choosing for the cost multiplier $p_0$ :

$$p_0 = -\frac{7}{9} \frac{\alpha^2 - 1}{\alpha^2} \frac{1}{f_{max}} < 0 \tag{50}$$

the minimal acceleration value $f_m$ is :

$$f_m = p_\Omega \dot{\Omega}_d f_{max} \tag{51}$$

The singular arc is then completely defined by the parameters $p_\Omega$ and $\dot{\Omega}_d$ and by the control Eqs. (44,45).

A singular arc may occur in the trajectory only if $0 < p_\Omega \dot{\Omega}_d < 1$ and it can be followed as long as the acceleration level Eq. (45) does not exceed $f_{max}$.

Similarly in the case $p_\Omega \dot{\Omega}_d < 0$, the sens of variation in Table 1 is reversed on each interval. The singular arc may occur only in an inclination range near 0 or 180 deg. With the adequate choice of $p_0$, the minimal acceleration level $f_0$ is :

$$f_0 = -p_\Omega \dot{\Omega}_d f_{max} \quad \text{with} \quad p_0 = -\frac{\alpha^2}{\alpha^2 - 7} \frac{1}{f_{max}} < 0 \tag{52}$$

A numerical procedure is necessary to detect the possible junctions between regular and singular arcs, and to find the corresponding values of the parameters $p_\Omega$ and $\dot{\Omega}_d$. In the case of a single singular arc spanning from $t_0$ to $t_f$, additional properties can be established.

Case of a single singular arc

From Property 3 the precession rate is constant along the singular arc so that : $\dot{\Omega}(t_0) = \dot{\Omega}(t_f) = \dot{\Omega}(t) = \dot{\Omega}_d$.

$\dot{\Omega}_d$ is thus equal to $\dot{\Omega}(t_0)$ which is known. There remains as only parameter $p_\Omega$.

The precession rate being constant throughout the transfer, the problem becomes insensitive to the initial and the final date since the relative configuration of the orbits does not evolve with the time. This is in accordance with Property 5 yielding that the cost total derivative Eqs. (18,20) nullify at the endpoints of a singular arc.

The cost of the singular solution can be assessed analytically starting from the velocity equation Eq. (4).

$$\dot{V} = -f \cos\beta \quad \Rightarrow \quad fdt = -\frac{dV}{\cos\beta} \tag{53}$$

Using $\dot{\Omega} = -kV^7 \cos I = \dot{\Omega}_d$ to eliminate V and Eq. (38) to eliminate β, we get a quadrature for the cost :

$$J = \int_{t_0}^{t_f} fdt = \pm \int_{I_0}^{I_f} \frac{1}{7} \left(\frac{\dot{\Omega}_d}{k \cos I}\right)^{\frac{1}{7}} \sqrt{\tan^2 I + \alpha^2} \, dI \tag{54}$$

The total cost of the singular solution can be assessed directly since it depends only on the initial and final inclinations, and on the initial precession rate. In particular it does not depend on the evolution of the acceleration level given by Eq. (45).

In order to define completely the singular solution, the value of $p_\Omega$ must be found so that the final conditions ($V_f$, $I_f$) are met at the prescribed final date $t_f$ when applying the control Eq. (44,45). Since the precession rate is constant along the trajectory : $\dot{\Omega}(t) = \dot{\Omega}(t_0)$, the velocity depends explicitly on the inclination. The state propagation Eq. (4) can be restricted to the inclination. The problem reduces to finding the unknown $p_\Omega$ such that the targeted inclination $I_f$ is reached at $t_f$.

This preliminary analysis shows that single singular solutions may exist. There remains open questions regarding the optimality of such solutions and also the existence of mixed solutions with possible junctions between regular and singular arcs.

## 2.6 Solution Recap

The Table 2 recaps the possible solutions of OCP Eq. (4), either purely regular, or purely singular or mixed.

|  | Abnormal | Normal | | |
| --- | --- | --- | --- | --- |
|  |  | $p_\Omega \neq 0$ | $p_\Omega = 0$ and $H \neq 0$ | $p_\Omega = 0$ and $H = 0$ |
| Regular | Edelbaum's solution $f = f_{max}$ | 3 sequences $f_{max} - 0 - f_{max}$ | Edelbaum's solution $f = f_{max}$ | Coast trajectory $f = 0$ |
| Singular | None | Explicit control $f, \beta(I)$ Cost independent on f Unknown $p_\Omega$ | None | Explicit state $V, I(\beta)$ Cost independent on f Unknown $f(t)$ |
| Mixed | None | To be investigated | None | None |

Table 2 : Possible solutions of OCP

# 3. Solution Method

This section deals with the numerical solution of the optimal control problem in the regular case. A shooting method is applied to the TPBVP issued from the PMP necessary conditions. The shooting method is initiated with a costate guess derived from a simplified problem.

## 3.1 Shooting Problem

The extremal analysis has shown (Property 6) that an optimal transfer in the regular case could be sought with 3 sequences and a control structure : $f_{max} - 0 - f_{max}$. The switching dates are denoted $t_1$ and $t_2$. The trajectory is obtained by integration of the state and costate Equations (4,6) from the initial date $t_0$ to the final date $t_f$. The acceleration level is set successively to $f_{max}$ from $t_0$ to $t_1$, to zero from $t_1$ to $t_2$ and to $f_{max}$ from $t_2$ to $t_f$. The thrust direction β is defined by Eq. (11).

The problem unknowns are the initial costate components and the switching dates. The final state is prescribed and the switching function must vanish at the switching dates. This shooting problem formulates as follows.

Shooting problem

Find the 3 initial costate components $p_V(t_0)$, $p_I(t_0)$, $p_\Omega(t_0)$ and the 2 switching dates $t_1$, $t_2$ in order to nullify the shooting function $F : R^5 \rightarrow R^5$ defined as :

$$F(p_V(t_0), p_I(t_0), p_\Omega(t_0), t_1, t_2) = (V(t_f) - V_f, I(t_f) - I_f, \Omega(t_f) - \Omega_f, S(t_1), S(t_2)) \qquad (55)$$

The shooting problem could be formulated with only 3 unknowns (initial costate) and 3 equations (final state), since the acceleration level can be determined from the switching function sign within the trajectory propagation. For the numerical resolution, it proves more efficient to add the 2 switching dates and the 2 switching conditions to the nonlinear systems.

The shooting method consists in solving the nonlinear system Eq. (55) by a Newton method. The main issue lies in the initial costate guess that must be sufficiently close to the solution in order to achieve the convergence of the Newton method. For that purpose the optimal control problem is approximated by a simplified nonlinear programming (NLP) problem that can be easily solved. The costate guess is then derived from the cost sensitivities of this NLP problem.

## 3.2 Split Edelbaum Strategy

The control structure of a regular solution is $f_{max} - 0 - f_{max}$. The trajectory is split into a first propelled transfer from the initial orbit ($V_0$,$I_0$) towards the drift orbit ($V_d$,$I_d$), a coast phase on the drift orbit and a second propelled transfer from the drift orbit ($V_d$,$I_d$) towards the final orbit ($V_f$,$I_f$). The fuel consumption is directly proportional to the durations of the two propelled transfers.

The simplifying approach consists in using the Edelbaum's analytical minimum-time solution to assess the two propelled transfers. The problem reduces to the drift orbit parameters ($V_d$,$I_d$) that should be sought in order to minimize the cost whilst meeting the final RAAN constraint. This simplified transfer strategy is sub-optimal wrt OCP, since the Edelbaum's solution does not account for the RAAN constraint. It should nevertheless be quite close to the OCP solution as long as the propelled duration remains small compared to the drift duration, meaning that the RAAN change is essentially performed during the drift phase. This assumption is the underlying motivation for the formulation of the OCP using the Edelbaum's dynamics.

The transfer strategy is depicted on the Figure 4, with the spiraling trajectories associated to the initial and final Edelbaum's transfers.

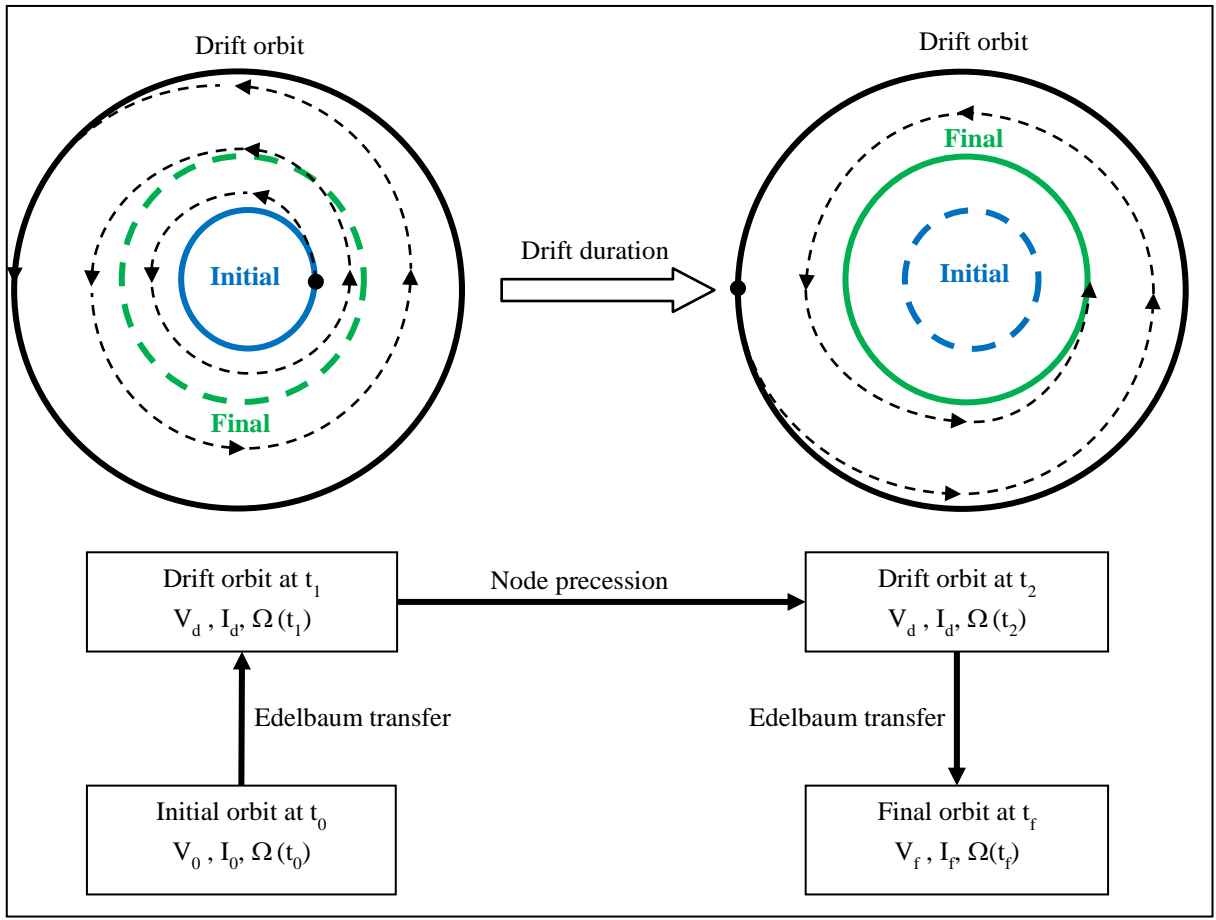

Figure 4 : Split Edelbaum Strategy

The simplified NLP problem corresponding to the Split Edelbaum Strategy (SES) formulates as follows.

SES problem

$$\min_{V_d, I_d} \Delta V \quad \text{s.t.} \quad \Omega(t_f) = \Omega_f \tag{56}$$

For both propelled transfers respectively from $(V_0, I_0)$ to $(V_d, I_d)$ and from $(V_d, I_d)$ to $(V_f, I_f)$, the Edelbaum's formulae [4,8] yield explicit assessments of the minimum duration, of the corresponding velocity impulse and also of the velocity V and inclination I evolutions wrt the time. The subscript "e" in the sequel refers to the Edelbaum's analytical solution. These formulae are used to assess on the one hand the transfer cost, on the other hand the final RAAN value as explained below.

The transfer cost is the sum of velocity impulses of the two propelled transfers respectively $\Delta V_{e1}$ and $\Delta V_{e2}$.

$$\Delta V = \Delta V_{e1} + \Delta V_{e2} \tag{57}$$

The Edelbaum's model considers a thrusting strategy that modifies the inclination at each revolution whilst leaving the RAAN unchanged. The RAAN evolution throughout the transfer is only due to the $J_2$ precession effect. The final RAAN value is obtained by integration of Eq. (1).

During the propelled transfers, the velocity V and the inclination I evolutions are given analytically by the Edelbaum's solution. These evolutions along both propelled transfers are denoted respectively $(V_{e1}, I_{e1})$ and $(V_{e2}, I_{e2})$. The RAAN evolution is assessed by quadrature for the propelled transfers, and explicitly for the drift phase since the precession rate is constant.

$$\Omega(t_1) = \Omega(t_0) + \int_{t_0}^{t_1} -kV_{e1}^7 \cos I_{e1} dt$$

$$\Omega(t_2) = \Omega(t_1) - kV_d^7 \cos I_d (t_2 - t_1) \tag{58}$$

$$\Omega(t_f) = \Omega(t_2) + \int_{t_2}^{t_f} -kV_{e2}^7 \cos I_{e2} dt$$

The cost for the SES problem Eq. (56) is thus assessed analytically whereas the constraint is assessed by 2 quadratures. No integration of differential equations is required.

### 3.3 SES Solution

The SES problem Eq. (56) with 2 unknowns and 1 equality constraint must be solved numerically using a nonlinear programming software. It can be initiated with a guess based on Edelbaum's analytical costates.

For each propelled transfer the endpoint costates provide indeed the cost derivatives wrt to the endpoint states. The drift orbit parameters represent the final state for the first propelled transfer (date $t_1$) and the initial state for the second propelled transfer (date $t_2$). Denoting respectively $(p_{eV1}, p_{eI1})$ and $(p_{eV2}, p_{eI2})$ the Edelbaum's costates of the two propelled transfers, we have :

$$\begin{cases} \dfrac{\partial(\Delta V_{e1})}{\partial V_d} = p_{eV1}(t_1) \\ \dfrac{\partial(\Delta V_{e1})}{\partial I_d} = p_{eI1}(t_1) \end{cases} \text{and} \quad \begin{cases} \dfrac{\partial(\Delta V_{e2})}{\partial V_d} = p_{eV2}(t_2) \\ \dfrac{\partial(\Delta V_{e2})}{\partial I_d} = p_{eI2}(t_2) \end{cases} \tag{59}$$

The gradient of the SES cost function Eq. (57) is thus :

$$\begin{cases} \dfrac{\partial \Delta V}{\partial V_d} = p_{eV1}(t_1) + p_{eV2}(t_2) \\ \dfrac{\partial \Delta V}{\partial I_d} = p_{eI1}(t_1) + p_{eI2}(t_2) \end{cases} \tag{60}$$

Assuming that the RAAN change is essentially achieved during the drift phase (which is the underlying motivation of the formulation Eq. (4)), the quadratures can be neglected in Eq. (58) and the RAAN constraint becomes analytical :

$$\dot{\Omega}_d = -kV_d^7 \cos I_d = \frac{\Omega_f - \Omega_0}{t_f - t_0} = C^{te} \tag{61}$$

Eq. (61) can be used to eliminate either $I_d$ or $V_d$. Choosing for instance $V_d$ as free parameter the cost gradient is :

$$\frac{d\Delta V}{dV_d} = \frac{\partial \Delta V}{\partial V_d} + \frac{\partial \Delta V}{\partial I_d}\frac{dI_d}{dV_d} = [p_{eV1}(t_1) + p_{eV2}(t_2)] + [p_{eI1}(t_1) + p_{eI2}(t_2)]\frac{7}{V_d \tan I_d} \tag{62}$$

The assessment of this gradient for any value of $V_d$ is completely analytical using first Eq. (61) to get $I_d$, and then the Edelbaum's formulae to get the endpoint costates $p_{eV1}$, $p_{eI1}$, $p_{eV2}$ and $p_{eI2}$. The optimal value $V_d$ must nullify the gradient Eq. (62). The problem reduces to a one unknown nonlinear equation which can be easily solved. The result is used as starting point for solving the SES problem Eq. (56).

It can be noticed that the difference $\Omega_f - \Omega_0$ in the second member of Eq. (61) is defined modulo $2\pi$. The SES problem Eq. (56) and similarly the OCP Eq. (4) have therefore several local minima corresponding to different drift precession rates $\dot{\Omega}_d$. In most cases the best solution can be guessed a priori depending on the relative configuration of the initial and final orbits in terms of RAAN values and of precession rates. Nevertheless the existence of local minima must be kept in mind when solving a particular instance of OCP Eq. (4) and in some cases different solutions must be compared in order to ensure that the global minimum is found.

### 3.4 OCP costate guess

For an optimal control problem with state vector X, costate vector P, cost function J and Hamiltonian function H, the following relationships hold [11].

$$P(t_0) = \frac{\partial J^*}{\partial X(t_0)} \quad , \quad H(t_0) = \frac{\partial J^*}{\partial t_0} \tag{63}$$

$J^*$ is the cost value for the extremal starting from $X(t_0)$ at the date $t_0$.

The initial costate is the partial derivative of the optimal cost wrt the initial state, whereas the initial Hamiltonian is the partial derivative of the optimal cost wrt the initial date. These properties are used to build a guess for the shooting problem Eq. (55). For that purpose the SES problem Eq. (56) is first solved for the reference initial conditions yielding the reference cost value. It is then solved again by varying the initial state components one by one and assessing the corresponding costate components by finite differences. An estimate of the Hamiltonian value is also assessed by varying the initial or the final date.

## 4. Application Case

The solution method is illustrated on a typical transfer of a debris removal mission. The initial and final orbits are circular and near sun-synchronous (from 800km/98deg to 900km/99deg).

### 4.1 Example Data

The Table 3 provides the initial and final conditions of the transfer. The transfer duration is bounded to 100 days. The maximal acceleration level is set to $3.5 \cdot 10^{-3}$ m/s$^2$.

|  | Initial orbit | Final orbit |
|---|---|---|
| Altitude (km) | 800.0 | 900.0 |
| Velocity (m/s) | 7450.0 | 7398.6 |
| Inclination (deg) | 98.00 | 99.00 |
| Precession rate (deg/day) | 0.917 | 0.982 |
| Initial RAAN (deg) at $t_0 = 0$ days | 0.0 | 30.0 |
| Final RAAN (deg) at $t_0 = 100$ days | 91.7 | 128.2 |

Table 3 : Initial and final orbits

The targeted orbit has a RAAN value of 30 deg at the transfer beginning with a precession rate of 0.982 deg/day. The RAAN value that must be targeted at the transfer end is then 128.2 deg ( = 30 + 0.982×100).

The RAAN gap of 30 deg must be caught up within 100 days. It is therefore expected that the drift orbit precession rate will be approximatively 0.3 deg/day higher than the target precession rate.

### 4.2 SES Solution

The SES problem Eq. (56) is solved using a reduced gradient optimizer (Airbus DS internal software). The initial values of the variables $(V_d, I_d)$ are obtained by solving the one variable equation (62). The following solution is obtained for the SES problem (in parenthesis the initial values) :

- Drift orbit velocity :     $V_d = 7664.0$ m/s     (7564.0 m/s)
- Drift orbit altitude :     $Z_d = 404.7$ km     (585.2 km)
- Drift orbit inclination :     $I_d = 99.20$ deg     (100.08 deg)
- Total velocity impulse :     $\Delta V = 598.1$ m/s
- Drift phase starting date :     $t_1 = 1.075$ days
- Drift phase ending date :     $t_2 = 99.137$ days

The transfer sequences are detailed in the Table 4.

|  | Date (day) | Altitude (km) | Velocity (m/s) | Inclination (deg) | RAAN (deg) | Precession (deg/day) | Impulse (m/s) |
|---|---|---|---|---|---|---|---|
| $t_0$ | 0.000 | 800.0 | 7450.0 | 98.00 | 0.00 | 0.917 | 0.0 |
| $t_1$ | 1.075 | 404.7 | 7664.0 | 99.20 | 1.18 | 1.284 | 328.7 |
| $t_2$ | 99.137 | 404.7 | 7664.0 | 99.20 | 127.26 | 1.284 | 328.7 |
| $t_f$ | 100.000 | 900.0 | 7398.6 | 99.00 | 128.20 | 0.982 | 598.1 |

Table 4 : SES Sequences

The drift precession rate is close to the expected value : 1.284 deg/day (≈ 0.982 + 0.3). This precession rate is obtained by an optimized combination of altitude decrease and inclination increase.

### 4.3 Cost Sensitivities

The SES problem is solved again with varying initial conditions in order to assess the cost sensitivities by centered finite differences. The results are presented in the Table 5 and converted in radians for the angles (inclination, RAAN). The initial velocity variation corresponds to an initial altitude variation of ±50 km. The sensitivity to the final date is also assessed.

|  | Variation | Cost ΔV (m/s) | Derivative (deg$^{-1}$) | Derivative (rad$^{-1}$) |
|---|---|---|---|---|
| Reference | 0 | 598.1 |  |  |
| Initial velocity (m/s) (Initial altitude km) | −25.9 (+50 km) | 615.5 | −0.630 (+0.327) | −0.630 (+0.327) |
|  | +26.0 (−50 km) | 582.8 |  |  |
| Initial inclination (deg) | +0.1 | 583.2 | −154.5 | −8852.2 |
|  | −0.1 | 614.1 |  |  |
| Initial RAAN (deg) | +5 | 539.3 | -13.27 | −760.31 |
|  | −5 | 672.0 |  |  |
| Final date (day) | +5 | 582.3 | -3.695 | -3.695 |
|  | −5 | 619.3 |  |  |

Table 5 : Cost derivatives

The cost partial derivatives wrt the initial state provide estimates of the initial costate components (the units are m/s , rad).

$$\begin{cases} p_V(t_0) \approx -0.630 \\ p_I(t_0) \approx -8852.2 \\ p_\Omega(t_0) \approx -760.31 \end{cases} \quad (64)$$

The Hamiltonian is estimated from the cost partial derivative wrt the final date (the units are m/s , day).

$$H = -\frac{\partial J^*}{\partial t_f} = -\frac{dJ^*}{dt_f} + \frac{\partial J^*}{\partial \Omega(t_f)} \times \frac{d\Omega(t_f)}{dt_f} = -\frac{dJ^*}{dt_f} + p_\Omega \dot{\Omega}_f = 3.695 + (-13.27) \times 0.982 = -9.336 \quad (65)$$

### 4.4 OCP Solution

The shooting problem Eq. (55) is solved by a Newton method (software Hybrd.f) initiated with the SES guess Eq. (64). This guess is already close to a zero of the shooting function as can be observed on the constraint values in the Table 6.

|  | Initial guess (m/s, rad, day) | Constraint value |
|---|---|---|
| Velocity | $p_V(t_0)$ = −0.630 | $V(t_f) - V_f$ = 5.7 m/s |
| Inclination | $p_I(t_0)$ = −8852.2 | $I(t_f) - I_f$ = −0.02 deg |
| RAAN | $p_\Omega(t_0)$ = −760.31 | $\Omega(t_f) - \Omega_f$ = −1.59 deg |
| Switching 1 | $t_1$ = 1.075 | $S(t_1)$ = 0.032 |
| Switching 2 | $t_2$ = 99.137 | $S(t_2)$ = 0.119 |

Table 6 : Shooting problem initial guess

The Newton method converges accurately in a few iterations yielding the shooting problem solution (Table 7).

|  | Initial guess (m/s, rad, day) | Constraint value |
|---|---|---|
| Velocity | $p_V(t_0)$ = −0.644 | $V(t_f) - V_f$ = $10^{-5}$ m/s |
| Inclination | $p_I(t_0)$ = −9215.9 | $I(t_f) - I_f$ = $10^{-8}$ deg |
| RAAN | $p_\Omega(t_0)$ = −816.97 | $\Omega(t_f) - \Omega_f$ = $10^{-6}$ deg |
| Switching 1 | $t_1$ = 1.092 | $S(t_1)$ = $10^{-8}$ |
| Switching 2 | $t_2$ = 99.114 | $S(t_2)$ = $10^{-8}$ |

Table 7 : Shooting problem solution

The sequences of the optimal transfer are detailed in the Table 8.

|  | Date (day) | Altitude (km) | Velocity (m/s) | Inclination (deg) | RAAN (deg) | Precession (deg/day) | Impulse (m/s) |
|---|---|---|---|---|---|---|---|
| $t_0$ | 0.000 | 800.0 | 7450.0 | 98.00 | 0.00 | 0.917 | 0,00 |
| $t_1$ | 1.092 | 407.1 | 7664.6 | 99.22 | 1.19 | 1.287 | 330.2 |
| $t_2$ | 99.114 | 407.1 | 7664.6 | 99.22 | 127.20 | 1.287 | 330.2 |
| $t_f$ | 100.000 | 900.0 | 7398.6 | 99.00 | 128.20 | 0.982 | 598.1 |

Table 8 : Optimal transfer sequences

In this exemple, the OCP solution is very close to the SES solution. This should be the case provided that a sufficient duration is allocated to the mission wrt the engine thrust level, so that the propelled durations remain small wrt the drift duration. If the duration allocated to the transfer becomes small, the Edelbaum's thrusting strategy (constant out of plane angle β with sign thrust at the antinodes) is no longer adapted to the mission. In such cases the optimal control problem must be addressed directly without averaging of the dynamics equations and the optimal thrust law should aim at changing simultaneously the inclination and the RAAN.

## 5. Conclusion

The minimum-fuel low-thrust transfer problem between circular orbit has been investigated considering a transfer strategy based on Edelbaum's averaged dynamics equations. The thrust is used to change the altitude and the inclination, whereas the RAAN change is realized passively by the natural precession due to the first zonal term. This control strategy requires a sufficient mission duration in order to benefit from the natural precession. It is especially suited for debris removal mission that are currently under study.

Abnormal and singular solutions have been considered and analytic formulae have been derived for the acceleration level in the singular case. Ongoing research focuses on the optimality of singular solutions and on the possible junctions between regular and singular legs.

For practical applications, only regular solutions are envisioned. Such regular trajectories are composed of 3 sequences. The first propelled sequence brings the vehicle on a circular drift orbit, the second coast sequence achieves the required RAAN change and the third propelled sequence completes the transfer to reach the targeted final altitude and inclination.

A shooting method is applied to solve the optimal control problem. The costate guess is derived from a nonlinear problem that approximate the transfer using Edelbaum's minimum-time solution, this nonlinear problem being itself initiated with a quasi-analytical solution. The guess is quite close to the optimal solution provided that a sufficient duration is allocated to the mission. The solution method proves robust enough for an automatic procedure that is used for the design of multiple debris removal missions.

**References**

[1] Docherty, S. Y., Macdonald, M. "Analytical Sun-Synchronous Low-Thrust Orbit Maneuvers," *Journal of Guidance, Control, and Dynamics*, Vol. 35, No.2, March-April 2012, pp. 681-686
[2] Edelbaum, T. N., "Propulsion Requirements for Controllable Satellites," *ARS Journal*, Vol 31, Aug. 1961, pp. 1079-1089.
[3] Leitmann, G., *Optimization technique with applications to Aerospace System*, Academic Press New-York, 1962.
[4] Kechichian, J.A., "Reformulation of Edelbaum's Low-Thrust Transfer Problem Using Optimal Control Theory," *Journal of Guidance, Control, and Dynamics*, Vol. 20, No.5, September-October 1997.
[5] Casalino, L., Colasurdo, G., "Improved Edelbaum's Approach to Optimize LEO-GEO Low-Thrust Transfers," AIAA/AAS Astrodynamics Specialist Conference, 16-19 August 2004.
[6] Conway, B.A., *Spacecraft Trajectory Optimization*, Cambridge University Press, 2010.
[7] Kluever, C.A., "Using Edelbaum's Method to Compute Low-Thrust Transfers with Earth-Shadow Eclipses," *Journal of Guidance, Control, and Dynamics*, Vol. 34, No.1, January-February 2011, pp. 300-303
[8] Chobotov, V., *Orbital mechanics*, 3rd Edition, AIAA Education Series, 2002
[9] Cerf, M., "Multiple Space Debris Collecting Mission – Optimal Mission Planning ," *Journal of Optimization Theory and Applications* ,Springer InterScience, 2015 , (DOI) 10.1007/s10957-015-0705-0
[10] Pontryagin, L., Boltyanskii, V., Gramkrelidze, R., Mischenko, E., *The mathematical theory of optimal processes*, Wiley Interscience, 1962.
[11] Trélat, E., *Contrôle optimal – Théorie et Applications*, Vuibert, 2005.

**Acronyms**

| | |
|---|---|
| EOR | Earth Orbit Raising |
| LEO | Low Earth Orbit |
| GEO | Geostationary Orbit |
| SSO | Sun Synchronous Orbit |
| RAAN | Right Ascension of Ascending Node |
| OCP | Optimal Control Problem |
| PMP | Pontryaguin Maximum Principle |
| TPBVP | Two Point Boundary Value Problem |
| NLP | Nonlinear Programming |
| SES | Split Edelbaum Strategy |